\documentclass[12pt]{article}

\usepackage[geometry]{ifsym}
\usepackage{verbatim}
\usepackage{fancyhdr}
\usepackage{graphicx}
\usepackage{mathrsfs}
\usepackage{amssymb}
\usepackage{cite}
\usepackage{epsfig}
\usepackage{pst-poly}     
\usepackage{pst-all}      

\newtheorem{pps}{Proposition}[section]

\newtheorem{lem}{Lemma}[section]
\newtheorem{thm}{Theorem}[section]

\newenvironment{pf}[1][Proof]{\noindent\textbf{#1.} }{\hfill\rule{1mm}{2mm}}

\makeatletter \@addtoreset{equation}{section} \makeatother

\begin{document}

\title{The total bondage number of grid graphs \thanks {The work was supported by NNSF
of China (No. 11071233).}}
\author
{Fu-Tao Hu\quad You Lu\quad Jun-Ming Xu\footnote{Corresponding
author:
xujm@ustc.edu.cn}\ \\ \\
{\small Department of Mathematics}  \\
{\small University of Science and Technology of China}\\
{\small Hefei, Anhui, 230026, China} }
\date{}
\maketitle

\begin{quotation}

\textbf{Abstract}: The total domination number of a graph $G$
without isolated vertices is the minimum number of vertices that
dominate all vertices in $G$. The total bondage number $b_t(G)$ of
$G$ is the minimum number of edges whose removal enlarges the total
domination number. This paper considers grid graphs. An $(n,m)$-grid
graph $G_{n,m}$ is defined as the cartesian product of two paths
$P_n$ and $P_m$. This paper determines the exact values of
$b_t(G_{n,2})$ and $b_t(G_{n,3})$, and establishes some upper bounds
of $b_t(G_{n,4})$.



\end{quotation}

\section{Introduction}

For notation and graph-theoretical terminology not defined here we
follow \cite{x03}. Specifically, let $G=(V,E)$ be an undirected
graph without loops and multi-edges, where $V=V(G)$ is the
vertex-set and $E=E(G)$ is the edge-set, which is a subset of
$\{xy|\ xy$ is an unordered pair of $V \}$. A graph $G$ is {\it
nonempty} if $E(G)\ne \emptyset$. Two vertices $x$ and $y$ are {\it
adjacent} if $xy\in E(G)$. For a vertex $x$, we call the vertices
adjacent to it the {\it neighbors} of $x$. We use $P_n$ and $C_n$ to
denote a path and a cycle of order $n$ throughout this paper.

A subset $D\subseteq V(G)$ is called a {\it dominating set} of $G$
if every vertex not in $D$ has at least one neighbor in $D$. The
{\it domination number} of $G$, denoted by $\gamma(G)$, is the
minimum cardinality among all dominating sets.

The domination is so an important and classic conception that it has
become one of the most widely studied topics in graph theory, and
also is frequently used to study properties of interconnection
networks. The early results on this subject have been surveyed and
detailed in the two excellent domination books by Haynes,
Hedetniemi, and Slater~\cite{hhs97, hhs98}. In the recent decade, a
large number of research papers on domination as well as related
topics appear in many scientific journals because of their
applications in many fields such as networks and so on.

A dominating set $D$ of a graph $G$ without isolated vertices is
called to be {\it total} if every vertex in $G$ has at least one
neighbor in $D$. The minimum cardinality among all total dominating
sets is called the {\it total domination number} of $G$, denoted by
$\gamma_t(G)$. It is clear that $\gamma(G)\leqslant\gamma_t(G)
\leqslant 2\gamma(G)$ for any graph $G$ without isolated vertices.

The total domination in graphs was introduced by Cockayne {\it et
al.}~\cite{cdh80} in 1980. The total domination in graphs has been
extensively studied in the literature. In 2009, Henning~\cite{h09}
gave a survey of selected recent results on this topic.

In 1990, Fink {\it et al.}~\cite{fjkr90} introduced the bondage
number as a parameter for measuring the vulnerability of the
interconnection network under link failure. The minimum dominating
set of sites plays an important role in the network for it dominates
the whole network with the minimum cost. So we must consider whether
its function remains good under attack. Suppose that someone such as
a saboteur does not know which sites in the network take part in the
dominating role, but does know that the set of these special sites
corresponds to a minimum dominating set in the related graph. Then
how many links does he have to attack so that the cost can not
remains the same in order to dominate the whole network? That
minimum number of links is just the bondage number.

The {\it bondage number} $b(G)$ of a nonempty graph $G$ is the
minimum number of edges whose removal from $G$ results in a graph
with larger domination number than $\gamma(G)$. Since the domination
number of every spanning subgraph of a nonempty graph $G$ is at
least as great as $\gamma(G)$, the bondage number of a nonempty
graph is well defined. Many results on this topic are obtained in
the literature. The exact values of the bondage numbers for some
graphs are determined, for example, a complete graph, a path, a
cycle, a complete $t$-partite graph~\cite{fjkr90}, a
tree~\cite{hr92, t97, tv00, hjvw98}, for the cartesian product of
two cycles $C_4\times C_n$~\cite{ksk05} and $C_3\times
C_n$~\cite{syj06}, and for other graphs~\cite{hx06, hx07, hx08}.
Some upper bounds of the bondage numbers for graphs are established,
see, for example, \cite{dhtv98, fjkr90, hr94, hp08, hx07a, hx08,
ls03, t97} for general graphs, \cite{cd06, frv03, ky00} for planar
graphs.
In particular, very recently, Hu and Xu~\cite{hx10} have
showed that the problem determining bondage number for general
graphs is NP-hard.

Following Fink et al., Kulli and Patwari~\cite{kp91} proposed the
concept of the total bondage number for a graph. The total bondage
number $b_t(G)$ of a graph $G$ is the minimum number of edges whose
removal results in a graph with total domination number larger than
$\gamma_t (G)$. If $b_t(G)$ does not exist, for example a star graph
$K_{1,n}$, we define $b_t(G)=\infty$. Kulli and Patwari~\cite{kp91}
calculated the exact values of $b_t(G)$ for some standard graphs
such as a cycle $C_n$ and a path $P_n$ for $n\geqslant 4$, a
complete bipartite graph $K_{m,n}$ and a complete $K_n$. Sridharan
{\it et al.}~\cite{ses07} showed that for any positive integer $k$
there exists a tree $T$ with $b_t(T)=k$. These authors also
established the upper bounds of $b_t(G)$ for a graph $G$ in terms of
its order. As far as we know, no much research work on the total
bondage number was reported in the literature except for the
above-mentioned. However, Hu and Xu~\cite{hx10} also showed that the
problem determining total bondage number for general graphs is
NP-hard.

An $(n,m)$-grid graph $G_{n,m}$ is the cartesian product $P_n\times
P_m$ of two paths $P_n$ and $P_m$. In this paper, we consider
$b_t(G_{n,m})$. Since $G_{1,m}\cong P_m$, we assume $n\ge 2$ under
our discussion. In 2002, Gravier~\cite{g02} determined $\gamma_t
(G_{n,m})$ for any $m\in\{1,2,3,4\}$, based on which we obtain the
following results.
 $$
 \begin{array}{rl}
 & b_{t}(G_{n,2})=\left\{ \begin{array}{l}
 1 \ {\rm if}\ n\equiv 0\,({\rm mod}\, 3),  \\
2 \ {\rm if}\ n\equiv 2\,({\rm mod}\, 3),\\
3 \ {\rm if}\ n\equiv 1\,({\rm mod}\, 3);
\end{array}
 \right.\\
 & b_{t}(G_{n,3})=1; \  b_{t}(G_{6,4})=2, {\rm and}\\
 & b_{t}(G_{n,4})\left\{ \begin{array}{ll}
 =1\ & {\rm if}\ n\equiv 1\,({\rm mod}\, 5)\ {\rm and}\ n\ne 6;\\
 =2\ & {\rm if}\ n\equiv 4\,({\rm mod}\, 5);\\
 \leq 3\ & {\rm if}\ n\equiv 2\,({\rm mod}\, 5);\\
 \leq 4\ & {\rm if}\ n\equiv 0,3\,({\rm mod}\, 5).\\
 \end{array} \right.
 \end{array}
 $$

The proofs of these results are in Section 3, Section 4 and Section
5, respectively. In Section 2, we give two preliminary results to be
used in our proofs.


\section{Preliminary results}

\begin{figure}[ht]
\begin{center}
\begin{pspicture}(-5,-.5)(5,3.5)

\cnode(-4,0){3pt}{1}\rput(-3.7,0){\scriptsize 1}
\cnode(-4,1){3pt}{2}\rput(-3.7,1){\scriptsize 2}
\cnode(-4,2){3pt}{3}\rput(-3.7,2){\scriptsize 3}
\cnode(-4,3){3pt}{4}\rput(-3.7,3){\scriptsize 4} \ncline{1}{2}
\ncline{2}{3} \ncline{3}{4}

\cnode(-2,.5){3pt}{1'}\rput(-1.7,0.5){\scriptsize 1}
\cnode(-2,1.5){3pt}{2'}\rput(-1.7,1.5){\scriptsize 2}
\cnode(-2,2.5){3pt}{3'}\rput(-1.7,2.5){\scriptsize 3}
\ncline{1'}{2'} \ncline{2'}{3'}

\rput(-3,1.5){$\times$}  \rput(-1,1.5){=}

\cnode(0,0){3pt}{11}\rput(.4,-.3){\scriptsize $x_{11}$}
\cnode(0,1.5){3pt}{12}\rput(.4,1.2){\scriptsize $x_{12}$}
\cnode(0,3){3pt}{13}\rput(.4,2.7){\scriptsize $x_{13}$}
\cnode(1.5,0){3pt}{21}\rput(1.9,-.3){\scriptsize $x_{21}$}
\cnode(1.5,1.5){3pt}{22}\rput(1.9,1.2){\scriptsize $x_{22}$}
\cnode(1.5,3){3pt}{23}\rput(1.9,2.7){\scriptsize $x_{23}$}
\cnode(3,0){3pt}{31}\rput(3.4,-.3){\scriptsize $x_{31}$}
\cnode(3,1.5){3pt}{32}\rput(3.4,1.2){\scriptsize $x_{32}$}
\cnode(3,3){3pt}{33}\rput(3.4,2.7){\scriptsize $x_{33}$}
\cnode(4.5,0){3pt}{41}\rput(4.9,-.3){\scriptsize $x_{41}$}
\cnode(4.5,1.5){3pt}{42}\rput(4.9,1.2){\scriptsize $x_{42}$}
\cnode(4.5,3){3pt}{43}\rput(4.9,2.7){\scriptsize $x_{43}$}

\ncline{11}{12} \ncline{12}{13} \ncline{10}{20} \ncline{11}{21}
\ncline{12}{22} \ncline{13}{23} \ncline{21}{22} \ncline{22}{23}
\ncline{20}{30} \ncline{21}{31} \ncline{22}{32} \ncline{23}{33}
\ncline[linewidth=2pt]{31}{32} \ncline[linewidth=2pt]{32}{33}
\ncline{30}{40} \ncline[linewidth=2pt]{31}{41} \ncline{32}{42}
\ncline[linewidth=2pt]{33}{43} \ncline[linewidth=2pt]{41}{42}
\ncline[linewidth=2pt]{42}{43} \ncline{40}{50} \ncline{41}{51}
\ncline{42}{52} \ncline{43}{53}

\end{pspicture}
\caption{\label{f1}\footnotesize  A $(4,3)$-grid graph
$G_{4,3}=P_4\times P_3$}
\end{center}
\end{figure}
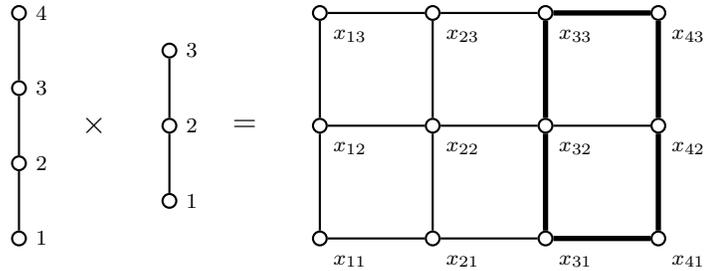

Throughout this paper, we assume that a path $P_n$ has the
vertex-set $V(P_n)=\{1,\cdots,n\}$. An $(n,m)$-grid graph $G_{n,m}$
is defined as the Cartesian product $G_{n,m}=P_n \times P_m$ with
vertex-set $V(G_{n,m})=\{x_{ij}|\ 1\leq i \leq n, 1\leq j \leq m\}$
and two vertices $x_{ij}$ and $x_{i'j\,'}$ being linked by an edge
if and only if either $i=i'\in V(P_n)$ and $jj\,'\in E(P_m)$, such
an edge is called a {\it vertical edge}, or $j=j\,'\in V(P_m)$ and
$ii'\in E(P_n)$, such an edge is called a {\it horizontal edge}. The
graph shown in Figure~\ref{f1} is a $(4,3)$-grid graph $G_{4,3}$. It
is clear, as a graphic operation, that the cartesian product
satisfies commutative associative law if identify isomorphic graphs,
that is, $G_{n,m}\cong G_{m,n}$.

The following notations continually appear in our proofs. For a
given integer $t$ with $t<n$, $G_{t,m}$ is a subgraph of $G_{n,m}$.
We use the notation $H_{n-t,m}$ to denote $G_{n,m}-G_{t,m}$, that
is, $H_{n-t,m}$ is a subgraph of $G_{n,m}$ induced by the set of
vertices $\{x_{ij}|\ t+1\leq i \leq n, 1\leq j \leq m\}$. Clearly,
$H_{n-t,m}\cong G_{n-t,m}$. For example, the graph shown in
Figure~\ref{f1} by heavy lines is a subgraph $H_{2,3}$ of $G_{4,3}$,
where $n=4, t=2$ and $m=3$.

Note that both $G_{0,m}$ and $H_{n-n,m}$ are nominal graphs. For
convenience of statements, we allow $G_{0,m}$ and $H_{n-n,m}$ to
appear in our proofs. If so, we specify their total dominating sets
are empty.

In Addition, let $Y_i=\{x_{ij}|\ 1\leq j \leq m\}$ for $1\leq i\leq
n$, called a set of {\it vertical vertices} in $G_{n,m}$.

We state some useful results on $\gamma_{t}(G_{n,m})$ to be used in
our proofs.

\begin{lem}\label{lem2.1}{\rm(Gravier~\cite{g02})}\
Let $n$ be a positive integer. Then

$\gamma_{t}(G_{n,2})=2\lfloor\frac{n+2}{3}\rfloor$ for $n\geq 1$;

$\gamma_{t}(G_{1,3})=2$ and $\gamma_{t}(G_{n,3})=n$ for $n\geq 2$;

$\gamma_t(G_{n,4})=\left\{
\begin{array}{l}
 \lfloor \frac{6n+8}{5}\rfloor \ \ \ \ \ \  {\rm if}\ n\equiv 1,2,4 \,({\rm mod}\, 5),  \\
\lfloor \frac{6n+8}{5}\rfloor+1 \ {\rm otherwise}
\end{array}
 \right. $ for $n\geq 4$.
 \end{lem}


\begin{lem}\label{lem2.2}
Let $D$ be a total dominating set of $G_{n,m}$. Then
$\gamma_{t}(G_{i,m})\leq |D\cap V(G_{i+1,m})|$ for $1\leq i \leq
n-1$ and $m\geq 2$.
\end{lem}

\begin{pf}
Let $D'=D\cap V(G_{i+1,m})$. If $i=1$, then the lemma holds clearly.
Assume $i\ge 2$ below.

If $D'\cap Y_{i+1}=\emptyset$, then $D'$ is a total dominating set
of $G_{i,m}$, and hence $\gamma_{t}(G_{i,m})\leq |D'|$.

Assume $D'\cap Y_{i+1}\ne\emptyset$ below. By the definition of a
total dominating set, it is possible that a vertex in $D'\cap Y_{i}$
is dominated only by vertices $D'\cap Y_{i+1}$. Let $A_i$ be the set
of vertices in $D'\cap Y_{i}$ dominated only by vertices in $D'\cap
Y_{i+1}$, and let $B_i=\{j|\ x_{ij}\in A_i\}$. Then
$D''=(D'\setminus Y_{i+1})\cup \{x_{(i-1)j}|\ j\in B_i\}$ is a total
dominating set of $G_{i,m}$ and $|D''|\leq |D'|$. Thus, we have
$\gamma_{t}(G_{i,m})\leq |D''|\leq |D'|$. The lemma follows.
\end{pf}


\begin{lem}\label{lem2.3}\textnormal{(Kulli and Patwari~\cite{kp91})}
For
a path $P_n$ with $n\geqslant 4$,
 $$
b_t(P_n)=\left\{\begin{array}{ll}
2&{\rm if}\ n\equiv 2\,({\rm mod}\,4);\\
1&{\rm otherwise}.
\end{array}\right.
 $$
\end{lem}

Since $G_{1,m}\cong P_m$ and $G_{n,1}\cong P_n$, by
Lemma~\ref{lem2.3}, we assume that if one of $n$ and $m$ is $1$,
then the other is at least $4$ when we consider the existence of
$b_t(G_{n,m})$.


\section{The total bondage number of $G_{n,2}$}

In this section, we determine the exact value of $b_t(G_{n,2})$ for
$n\geq 2$. Since the computation of $b_t(G_{n,2})$ strongly depends
on the value of $\gamma_t(G_{n,2})$ in Lemma~\ref{lem2.1}, the
process of our proofs consists of several lemmas according to the
value of $n$ modulo $3$.

\begin{lem}\label{lem3.1}
$b_{t}(G_{n,2})\leq 2$ for $n\not\equiv 1\,({\rm mod}\, 3)$.
\end{lem}

\begin{pf}
By Lemma \ref{lem2.1}, we have
 \begin{equation}\label{e3.1}
 \gamma_{t}(G_{n-1,2})=\gamma_{t}(G_{n,2})\ {\rm if}\
 n\not\equiv 1\,({\rm mod}\, 3).
 \end{equation}
Let $B=\{x_{(n-1)1}x_{n1},x_{(n-1)2}x_{n2}\}\subset E(G_{n,2})$, and
let $H=G_{n,2}-B$. By (\ref{e3.1}), we have
$$
\gamma_{t}(H)=2+\gamma_{t}(G_{n-1,2}) \geq 1+\gamma_{t}(G_{n,2}),
$$
which implies that $b_{t}(G_{n,2})\leq |B|=2$.
\end{pf}


\begin{lem}\label{lem3.2}
If $n\equiv 1\,({\rm mod}\, 3)$, then
$\gamma_{t}(G_{n,2}-x_{nj})=\gamma_{t}(G_{n,2})-1$ for each $j=1,2$.
\end{lem}

\begin{pf}
Without loss of generality, we only consider the case $j=1$. By the
hypothesis, $n\geq 4$. It can be direct check that the lemma holds
for $n=4$. Assume $n\geq 7$ below. Let $D$ be a minimum total
dominating set of $G_{n,2}-x_{n1}$. We need to show
$|D|=\gamma_{t}(G_{n,2})-1$.

We consider a subgraph $G_{n-4,2}$ and let $D'$ be a minimum total
dominating set of $G_{n-4,2}$.
By Lemma~\ref{lem2.1}, $|D'|=2\lfloor
\frac{n-4+2}{3}\rfloor=2\lfloor \frac{n-2}{3}\rfloor$. Clearly,
$D'\cup \{x_{(n-1)2},x_{(n-2)2},x_{(n-3)2}\}$ is a total dominating
set of $G_{n,2}-x_{n1}$. Since $n\equiv 1\,({\rm mod}\, 3)$, we have
$2\lfloor
 \frac{n+4}{3}\rfloor=2\lfloor
\frac{n+2}{3}\rfloor=\gamma_{t}(G_{n,2})$ by Lemma~\ref{lem2.1}. It
follows that
 $$
 \begin{array}{rl}
 |D|&\leq |D'|+3=2\lfloor
\frac{n-2}{3}\rfloor+3\\
 &=2\lfloor\frac{n+4}{3}\rfloor-1=2\lfloor
 \frac{n+2}{3}\rfloor-1\\
 &=\gamma_{t}(G_{n,2})-1,
 \end{array}
 $$
that is,
 \begin{equation}\label{e3.2}
|D|\leq \gamma_{t}(G_{n,2})-1.
 \end{equation}

We now prove that $|D|\geq \gamma_{t}(G_{n,2})-1$. If one of
$x_{(n-1)1}$ and $x_{n2}$ belongs to $D$, then $D$ is a total
dominating set of $G_{n,2}$. By (\ref{e3.2}), we can deduce a
contradiction as follows. $\gamma_{t}(G_{n,2})\leq |D|\leq
\gamma_{t}(G_{n,2})-1$. It follows that neither of $x_{(n-1)1}$ and
$x_{n2}$ belongs to $D$. Since $D$ is a total dominating set of
$G_{n,2}-x_{n1}$, the vertex $x_{(n-1)2}$ must be in $D$ to dominate
$x_{n2}$.
Thus $D\cup \{x_{(n-1)1}\}$ is a total dominating set of
$G_{n,2}$, and so
 $$
 |D|=|D\cup
\{x_{(n-1)1}\}|-1\geq\gamma_{t}(G_{n,2})-1.
$$

The lemma follows.
\end{pf}


\begin{lem}\label{lem3.3}
$b_{t}(G_{n,2})=1$ for $n\equiv 0\,({\rm mod}\, 3)$.
\end{lem}

\begin{pf}
We only need to show
 \begin{equation}\label{e3.3}
 \gamma_{t}(G_{n,2}-x_{(n-1)1}x_{n1})\geq \gamma_{t}(G_{n,2})+1.
 \end{equation}

Let $H=G_{n,2}-x_{(n-1)1}x_{n1}$ and $D$ be a minimum total
dominating set of $H$. Then the vertex $x_{n2}$ must be in $D$
otherwise $D$ can not dominate the vertex $x_{n1}$ in $H$. Moreover,
$D$ is either a total dominating set of $G_{n+1,2}$ if $x_{n1}$ is
in $D$ or a total dominating set of $G_{n+1,2}-x_{(n+1)1}$ if
$x_{n1}$ is not in $D$. Since $n\equiv 0\,({\rm mod}\, 3)$, we have
$n+1\equiv 1\,({\rm mod}\, 3)$ and
$\lfloor\frac{n}{3}\rfloor=\lfloor\frac{n+2}{3}\rfloor$. By Lemma
\ref{lem2.1} and Lemma \ref{lem3.2}, we have
 $$
 \begin{array}{rl}
 \gamma_{t}(H)=|D|& \geq \gamma_{t}(G_{n+1,2})-1
 =2\lfloor\frac{n+1+2}{3}\rfloor-1=2\lfloor\frac{n}{3}\rfloor+1\\
 &=2\lfloor\frac{n+2}{3}\rfloor+1=\gamma_{t}(G_{n,2})+1.
 \end{array}
 $$
The lemma follows.
\end{pf}


\begin{lem}\label{lem3.4}
$b_{t}(G_{n,2})=2$ for $n\equiv 2\,({\rm mod}\, 3)$.
\end{lem}

\begin{pf}
To prove the lemma, we only need to show $b_{t}(G_{n,2})\geq 2$ by
Lemma \ref{lem3.1}. To this end, we only need to show
$\gamma_{t}(G_{n,2}-e)=\gamma_{t}(G_{n,2})$ for any edge $e$ in
$G_{n,2}$. Let $e$ be any edge in $G_{n,2}$. We only need to prove
that $\gamma_{t}(G_{n,2}-e)\leq \gamma_{t}(G_{n,2})$ since
$\gamma_{t}(G_{n,2})\leq \gamma_{t}(G_{n,2}-e)$ clearly. We attain
this aim by constructing a total dominating set $D$ of $G_{n,2}-e$
such that $|D|=2\left\lfloor\frac{n+2}{3}\right\rfloor$, which means
$|D|=\gamma_{t}(G_{n,2})$ by Lemma~\ref{lem2.1}.

We consider two cases according as that $e$ is vertical or
horizontal.

Suppose that $e$ is a vertical edge $e=x_{i1}x_{i2}$, where $1\leq
i\leq n$. Let
 $$
 D=\left\{\begin{array}{ll}
 \{x_{kj}: k\equiv 1 \,({\rm mod}\, 3), j=1,2\}\ & {\rm if}\ i\equiv 2\,({\rm mod}\,3);\\
 \{x_{kj}: k\equiv 2\,({\rm mod}\, 3), j=1,2\}\ & {\rm otherwise}.
 \end{array}\right.
 $$
Then $D$ is a total dominating set of $G_{n,2}-e$ and
$|D|=2\left\lfloor\frac{n+2}{3}\right\rfloor$.

Suppose now that $e$ is a horizontal edge, maybe
$e=x_{i1}x_{(i+1)1}$ or $e=x_{i2}x_{(i+1)2}$, where $1\leq i\leq
n-1$. Without loss of generality, set $e=x_{i1}x_{(i+1)1}$. We
consider two subcases to construct $D$, respectively.

Assume $i\not\equiv 1\,({\rm mod}\, 3)$. Let
 $$
 D=\left\{\begin{array}{ll}
 \{x_{kj}: k\equiv 2 \,({\rm mod}\, 3), j=1,2\}\ & {\rm if}\ i\equiv 0\,({\rm mod}\, 3);\\
 \{x_{kj}: k\equiv 1\,({\rm mod}\, 3), j=1,2\}\ & {\rm if}\ i\equiv 2\,({\rm mod}\, 3).
 \end{array}\right.
 $$
Then $D$ is a total dominating set of $G_{n,2}-e$ and
$|D|=2\lfloor\frac{n+2}{3}\rfloor$.

\vspace{10pt}

\begin{center}
\begin{figure}[ht]
\begin{pspicture}(-0.5,-0.)(8.5,3)
\psset{unit=1.3cm}

\cnode(1,1){2pt}{11}\rput(1,0.8){\scriptsize $x_{11}$}
\cnode(1,2){2pt}{12}\rput(1,2.2){\scriptsize $x_{12}$}
\cnode(2,1){2pt}{21} \cnode(2,2){2pt}{22} \cnode(3,1){2pt}{31}
\cnode(3,2){2pt}{32} \cnode(4,1){2pt}{41}\rput(4,0.8){\scriptsize
$x_{(i-1)1}$} \cnode(4,2){2pt}{42}\rput(4,2.2){\scriptsize
$x_{(i-1)2}$} \cnode(5,1){2pt}{51}\rput(5,0.8){\scriptsize $x_{i1}$}
\cnode*(5,2){2pt}{52}\rput(5,2.2){\scriptsize $x_{i2}$}
\cnode(6,1){2pt}{61}\rput(6,0.8){\scriptsize $x_{(i+1)1}$}
\cnode*(6,2){2pt}{62}\rput(6,2.2){\scriptsize $x_{(i+1)2}$}
\cnode(7,1){2pt}{71}\rput(7,0.8){\scriptsize $x_{(i+2)1}$}
\cnode(7,2){2pt}{72}\rput(7,2.2){\scriptsize $x_{(i+2)2}$}
\cnode(8,1){2pt}{81} \cnode(8,2){2pt}{82} \cnode(9,1){2pt}{91}
\cnode(9,2){2pt}{92} \cnode(10,1){2pt}{101}\rput(10,0.8){\scriptsize
$x_{n1}$} \cnode(10,2){2pt}{102}\rput(10,2.2){\scriptsize $x_{n2}$}

\ncline{11}{12}  \ncline{11}{21} \ncline{12}{22} \ncline{21}{22}
\ncline[linestyle=dotted]{21}{31} \ncline[linestyle=dotted]{22}{32}
\ncline{31}{32} \ncline{31}{41} \ncline{32}{42} \ncline{41}{42}
\ncline{41}{51} \ncline{42}{52} \ncline{51}{52}
\ncline[linewidth=2pt]{51}{61}\rput(5.5,1.1){\scriptsize $e$}
\ncline{52}{62} \ncline{61}{62} \ncline{61}{71} \ncline{62}{72}
\ncline{71}{72} \ncline{71}{81} \ncline{72}{82} \ncline{81}{82}
\ncline[linestyle=dotted]{81}{91} \ncline[linestyle=dotted]{82}{92}
\ncline{91}{92} \ncline{91}{101} \ncline{92}{102} \ncline{101}{102}

\rput[b]{-90}(.8,0.5){$\left.\rule{0mm}{23mm}\right\}~\rput{90}{_{G_{i-1,2}}}$}
\rput[b]{-90}(6.9,0.5){$\left.\rule{0mm}{23mm}\right\}~\rput{90}{_{H_{n-(i+1),2}}}$}
\end{pspicture}
\caption{\label{f2}\footnotesize  Two subgraphs $G_{i-1,2}$ and
$H_{n-(i+1),2}$ of $G_{n,2}$}
\end{figure}
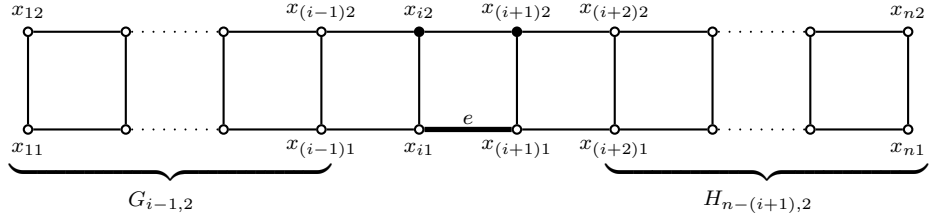
\end{center}

Assume now $i\equiv 1\,({\rm mod}\, 3)$. We consider $G_{i-1,2}$ and
$H_{n-(i+1),2}$ (see Figure~\ref{f2}). Let $D'$ and $D''$ be minimum
total dominating sets of $G_{i-1,2}$ and $H_{n-(i+1),2}$,
respectively. Then $D=D'\cup D''\cup \{x_{i2},x_{(i+1)2}\}$ is a
total dominating set of $G_{n,2}-e$. Note $D'=\emptyset$ if $i=1$
and $D''=\emptyset$ if $i=n-1$. Thus, by Lemma~\ref{lem2.1},
 $$
 \begin{array}{rl}
 |D|&=|D'|+|D''|+2\\
    &=2\left\lfloor\frac{i-1+2}{3}\right\rfloor+2
      \left\lfloor\frac{n-i-1+2}{3}\right\rfloor+2\\
    &=2\left\lfloor\frac{n+2}{3}\right\rfloor.
    \end{array}
 $$

The lemma follows.
\end{pf}


\begin{lem}\label{lem3.5}
$b_{t}(G_{n,2})=3$ for $n\equiv 1\,({\rm mod}\, 3)$.
\end{lem}

\begin{pf}
Since $n\equiv 1\,({\rm mod}\, 3)$, $n-1\equiv 0\,({\rm mod}\, 3)$.
By (\ref{e3.3}), for the edge $e_0=x_{(n-2)1}x_{(n-1)1}$, we have
 \begin{equation}\label{e3.4}
 \gamma_{t}(G_{n-1,2}-e_0)\geq \gamma_{t}(G_{n-1,2})+1.
 \end{equation}
Choose other two edges $e_1, e_2$ in $G_{n,2}$, where
$e_1=x_{(n-1)1}x_{n1}$ and $e_2=x_{(n-1)2}x_{n2}$. Let
$H=G_{n,2}-\{e_0, e_1, e_2\}$. Then
$H=(G_{n-1,2}-e_0)+H_{n-(n-1),2}$ and any total dominating set of
$H$ must contain vertices $x_{n1}$ and $x_{n2}$. By (\ref{e3.4}) and
Lemma~\ref{lem2.1}, we have
 $$
 \begin{array}{rl}
\gamma_{t}(H)
&=\gamma_{t}(G_{n-1,2}-e_0)+2\\
&\geq \gamma_{t}(G_{n-1,2})+1+2\\
&=2\left\lfloor\frac{n-1+2}{3}\right\rfloor+3\\
&=2\left\lfloor\frac{n+4}{3}\right\rfloor+1
=2\left\lfloor\frac{n+2}{3}\right\rfloor+1\\
&=\gamma_{t}(G_{n,2})+1,
\end{array}
$$
which implies $b_{t}(G_{n,2})\leq 3$.


Now we prove $b_{t}(G_{n,2})\geq 3$. To the end, let $e_1$ and $e_2$
be any two edges in $G_{n,2}$, and $H=G_{n,2}-\{e_1,e_2\}$. We only
need to prove $\gamma_{t}(H)\leq \gamma_{t}(G_{n,2})$. We consider
three cases.

\begin{description}
\item [Case 1] Both $e_1$ and $e_2$ are  vertical edges.

Let $e_1=x_{i1}x_{i2}$, $e_2=x_{j1}x_{j2}$, $i<j$, and let
 $$
 D=\left\{\begin{array}{ll}
 \{x_{kl}: k\equiv 1\,({\rm mod}\, 3), l=1,2\}\ & {\rm if}\ i,j\not\equiv 1\,({\rm mod}\, 3);\\
 \{x_{kl}: k\equiv 2\,({\rm mod}\, 3), l=1,2\}\cup\{x_{(n-1)1},x_{(n-1)2}\}\
 & {\rm if}\ i,j\not\equiv 2\,({\rm mod}\, 3);\\
 \{x_{kl}: k\equiv 0\,({\rm mod}\, 3), l=1,2\}\cup \{x_{21},x_{22}\} \ &
 {\rm otherwise}.
 \end{array}\right.
 $$
Then $D$ is a total dominating set of $H$ and $\gamma_{t}(H)\leq
|D|=2\left\lfloor\frac{n+2}{3}\right\rfloor$. By Lemma~\ref{lem2.1},
$|D|=\gamma_{t}(G_{n,2})$. Thus, for two vertical edges $e_1$ and
$e_2$, we have
 \begin{equation}\label{e3.5}
 \gamma_{t}(H)=\gamma_{t}(G_{n,2}-\{e_1,e_2\})\leq \gamma_{t}(G_{n,2}).
 \end{equation}

\item [Case 2] One of $e_1$ and $e_2$ is horizontal and the other is vertical.

Without loss of generality, suppose that $e_1$ is horizontal and
$e_2$ is vertical, and let $e_1=x_{i1}x_{(i+1)1}$ and
$e_2=x_{j1}x_{j2}$, $1\leq i\leq n-1$ and $1\leq j\leq n$. We will
prove $\gamma_{t}(H)\leq \gamma_{t}(G_{n,2})$.

Consider $G_{i,2}$ and $H_{n-i,2}$. Then both $G_{i,2}$ and
$H_{n-i,2}$ do not contain the edge $e_1$.
There are several subcases.

If $i\equiv 2\,({\rm mod}\, 3)$, then $n-i\equiv 2\,({\rm mod}\, 3)$
since $n\equiv 1\,({\rm mod}\, 3)$. By Lemma \ref{lem3.4},
$b_t(G_{i,2})=2=b_t(H_{n-i,2})$, which implies
$\gamma_{t}(G_{i,2})=\gamma_{t}(G_{i,2}-e_2)$ if $e_2$ is in
$G_{i,2}$, and $\gamma_{t}(H_{n-i,2})=\gamma_{t}(H_{n-i,2}-e_2)$ if
$e_2$ is in $H_{n-i,2}$. No matter which case arises, by
Lemma~\ref{lem2.1}, we have
 \begin{equation}\label{e3.6}
\begin{array}{rl}
 \gamma_{t}(H)&\leq \gamma_{t}(G_{i,2})+\gamma_{t}(H_{n-i,2})\\
 &=2\left\lfloor\frac{i+2}{3}\right\rfloor+2\left\lfloor\frac{n-i+2}{3}\right\rfloor\\
 &\leq 2\left\lfloor\frac{n+2}{3}\right\rfloor=\gamma_{t}(G_{n,2}).
 \end{array}
 \end{equation}

If $i\equiv 1\,({\rm mod}\, 3)$ and $j\le i$, then $e_2$ is in
$G_{i,2}$. By (\ref{e3.5}),
$\gamma_{t}(G_{i,2}-e_2)\leq\gamma_{t}(G_{i,2})$. Thus, the
inequalities (\ref{e3.6}) hold.

If $i\equiv 0\,({\rm mod}\, 3)$ and $j\geq i+1$, then $n-i\equiv
1\,({\rm mod}\, 3)$ and $e_2$ is in $H_{n-i,2}$. Since
$H_{n-i,2}\cong G_{n-i,2}$, by (\ref{e3.5}), we have
$\gamma_{t}(H_{n-i,2}-e_2)\leq\gamma_{t}(H_{n-i,2})$. Thus, the
inequalities (\ref{e3.6}) hold.

The remainder is the case either $i\equiv 1\,({\rm mod}\, 3)$ and
$j\ge i+1$ or $i\equiv 0\,({\rm mod}\, 3)$ and $j\le i$.
Essentially, the two cases are the same by replacing $n-i$ for $i$.
We only consider the latter case, that is, $i\equiv 0\,({\rm mod}\,
3)$ and $j\leq i$.

If $j=i$, let $D=\{x_{kl}: k\equiv 2\,({\rm mod}\, 3),
l=1,2\}\cup\{x_{n1},x_{n2}\}$, then $D$ is a total dominating set of
$H$, and so,
 $$
 \gamma_{t}(H)\leq |D|=
 2\left\lfloor\frac{n+2}{3}\right\rfloor=\gamma_{t}(G_{n,2}).
 $$

We now assume $j<i$. Consider $G_{i-1,2}$ and $H_{n-(i+1),2}$.
Let $D'$ be a minimum total dominating set of $G_{i-1,2}-e_2$, and
$D''$ be a minimum total dominating set of $H_{n-(i+1),2}$. Then
$D=D'\cup D''\cup \{x_{i2},x_{(i+1)2}\}$ is a total dominating set
of $H$. Since $i-1\equiv 2\,({\rm mod}\, 3)$,
$\gamma_{t}(G_{i-1,2})=\gamma_{t}(G_{i-1,2}-e_2)$ by Lemma
\ref{lem3.4}. $H_{n-(i+1),2}$ contains neither $e_1$ nor $e_2$. By
Lemma~\ref{lem2.1}, we have
 $$
 \begin{array}{rl}
 \gamma_{t}(H) &\leq |D|=|D'|+|D''|+2\\
 &\leq 2+\gamma_{t}(G_{i-1,2})+\gamma_{t}(H_{n-(i+1),2})\\
 &=2+2\lfloor\frac{i+1}{3}\rfloor+2\lfloor\frac{n-i+1}{3}\rfloor\\
 & =2\lfloor\frac{n+2}{3}\rfloor=\gamma_{t}(G_{n,2}).
 \end{array}
 $$

\item [Case 3] Both $e_1$ and $e_2$ are horizontal edges.

Without loss of generality, let $e_1=x_{i1}x_{(i+1)1}$ and
$e_2=x_{kj}x_{(k+1)j}$ are two distinct horizontal edges, where
$1\leq j\leq 2$ and $i\leq k<n$, and $j=2$ if $i=k$. To prove
$\gamma_{t}(H)\leq \gamma_{t}(G_{n,2})$, we consider three subcases.

\begin{description}

\item[Subcase 3.1] $k=i$.

In this subcase, $e_2=x_{i2}x_{(i+1)2}$, $H$ is disconnected and has
exact two connected components $G_{i,2}$ and $H_{n-i,2}$. Since both
$G_{i,2}$ and $H_{n-i,2}$ contain neither of $e_1$ and $e_2$, we
have
$\gamma_{t}(H)=\gamma_{t}(G_{i,2})+\gamma_{t}(H_{n-i,2})=\gamma_{t}(G_{n,2})$ by Lemma~\ref{lem2.1}.

\item[Subcase 3.2] $k=i+1$.

In this subcase, $G_{i,2}$ and $H_{n-i-1,2}$ contain neither $e_1$
nor $e_2$.

If $i\equiv 0$ or $1\,({\rm mod}\, 3)$, let $D'$ be a minimum total
dominating set of $G_{i-1,2}$, and $D''$ be a minimum total
dominating set of $H_{n-i-1,2}$, then $D=D'\cup D''\cup
\{x_{i2},x_{(i+1)2}\}$ is a total dominating set of $H$. Note
$D'=\emptyset$ if $i=1$. By Lemma~\ref{lem2.1}, we have
 $$
 \begin{array}{rl}
 \gamma_{t}(H)& \leq |D|\leq 2+\gamma_{t}(G_{i-1,2})+\gamma_{t}(H_{n-i-1,2})\\
 &=2+2\lfloor\frac{i+1}{3}\rfloor+2\lfloor\frac{n-i+1}{3}\rfloor\\
 &=2\lfloor\frac{n+2}{3}\rfloor=\gamma_{t}(G_{n,2}).
 \end{array}
 $$

If $i\equiv 2\,({\rm mod}\, 3)$, let $D'$ be a minimum total
dominating set of $G_{i-2,2}$, and $D''$ be a minimum total
dominating set of $H_{n-i-2,2}$, then $D=D'\cup D''\cup
\{x_{(i-1)p},x_{ip},x_{(i+1)p},x_{(i+2)p}\}$, where $p=3-j$, is a
total dominating set of $H$. Note $D'=\emptyset$ if $i=2$ and
$D''=\emptyset$ if $i=n-2$. By Lemma~\ref{lem2.1}, we have
 $$
 \begin{array}{rl}
 \gamma_{t}(H)& \leq |D|\leq 4+\gamma_{t}(G_{i-2,2})+\gamma_{t}(H_{n-i-2,2})\\
 &=4+2\lfloor\frac{i}{3}\rfloor+2\lfloor\frac{n-i}{3}\rfloor\\
 & =2\lfloor\frac{n+2}{3}\rfloor=\gamma_{t}(G_{n,2}).
 \end{array}
 $$

\item[Subcase 3.3] $k>i+1$.

In this case, $e_2$ is in $H_{n-i,2}$.

If $i\equiv 2\,({\rm mod}\, 3)$, then $n-i\equiv 2\,({\rm mod}\, 3)$
as $n\equiv 1\,({\rm mod}\, 3)$. Thus, $b_t(H_{n-i,2})=2$ by Lemma
\ref{lem3.4}, which implies
$\gamma_{t}(H_{n-i,2}-e_2)=\gamma_{t}(H_{n-i,2})$. If $i\equiv
0\,({\rm mod}\, 3)$, then $n-i\equiv 1\,({\rm mod}\, 3)$. Since
$e_2$ is a horizontal edge in $H_{n-i,2}$, by Subcase 3.1,
we also have $\gamma_{t}(H_{n-i,2}-e_2)=\gamma_{t}(H_{n-i,2})$.
Thus, when $i\not\equiv 1\,({\rm mod}\, 3)$, we have
 $$
 \begin{array}{rl}
 \gamma_{t}(H)&\leq \gamma_{t}(G_{i,2})+\gamma_{t}(H_{n-i,2}-e_2)\\
 & =\gamma_{t}(G_{i,2})+\gamma_{t}(H_{n-i,2})\\
 & =\gamma_{t}(G_{n,2})
 \end{array}
 $$

If $k\not\equiv 0\,({\rm mod}\, 3)$ then, by replacing $G_{k,2}$ and
$H_{n-k,2}$ by $H_{n-i,2}$ and $G_{i,2}$, respectively, we still
have
 $$
 \begin{array}{rl}
 \gamma_{t}(H)&\leq \gamma_{t}(G_{k,2}-e_1)+\gamma_{t}(H_{n-k,2})\\
 & =\gamma_{t}(G_{k,2})+\gamma_{t}(H_{n-k,2})\\
 &=\gamma_{t}(G_{n,2}).
 \end{array}
 $$

Now, we assume $i\equiv 1\,({\rm mod}\, 3)$ and $k\equiv 0\,({\rm
mod}\, 3)$. We consider three subgraphs $G_{i-1,2}$, $H_{n-(k+1),2}$
and $H_{k-1-(i+1),2}$.

Let $D'$ be a minimum total dominating set of $G_{i-1,2}$, $D''$ be
a minimum total dominating set of $H_{n-(k+1),2}$, and $D'''$ be a
minimum total dominating set of $H_{k-1-(i+1),2}$. Then $D=D'\cup
D''\cup D'''\cup \{x_{i2},x_{(i+1)2},x_{kp},x_{(k+1)p}\}$, where
$p=3-j$, is a total dominating set of $H$. Note $D'=\emptyset$ if
$i=1$, $D''=\emptyset$ if $k=n-1$ and $D'''=\emptyset$ if $k=i+2$.
By Lemma~\ref{lem2.1}, we have
 $$
 \begin{array}{rl}
 \gamma_{t}(H) & \leq |D|\leq 4+\gamma_{t}(G_{i-1,2})+\gamma_{t}(H_{n-k-1,2})+ \gamma_{t}(H_{k-i-2,2})\\
 & =4+2\lfloor\frac{i+1}{3}\rfloor+2\lfloor\frac{n-k-1+2}{3}\rfloor+2\lfloor\frac{k-i-2+2}{3}\rfloor\\
 &=4+2\lfloor\frac{n-4}{3}\rfloor=2\lfloor\frac{n+2}{3}\rfloor\\
 &=\gamma_{t}(G_{n,2}).
 \end{array}
 $$

\end{description}

\end{description}

Summing up all cases, we prove the lemma.
\end{pf}

\vskip6pt

According to the above lemmas, we can state our results in this
section as follows.

\begin{thm}\label{thm3.1} For any integer $n\geq 2$,
 $$
 b_{t}(G_{n,2})=\left\{ \begin{array}{l}
 1 \ {\rm if}\ n\equiv 0\,({\rm mod}\, 3),  \\
2 \ {\rm if}\ n\equiv 2\,({\rm mod}\, 3),\\
3 \ {\rm if}\ n\equiv 1\,({\rm mod}\, 3).
\end{array}
 \right.
 $$
\end{thm}


\section{The total bondage number of $G_{n,3}$}

In this section, we will determine $b_t(G_{n,3})=1$ for $n\geq 2$.
In this case, $Y_i=\{x_{ij}|\ 1\leq j \leq 3\}$ for $1\leq i\leq n$.

\begin{lem}\label{lem4.1}
Let $D$ be a minimum total dominating set of $G_{n,3}$. Then $|D\cap
Y_1|\leq 2$ and $|D\cap Y_n|\leq 2$.
\end{lem}

\begin{pf}
Without loss of generality, we only show $|D\cap Y_n|\leq 2$. By
contradiction, suppose that there exists a minimum total dominating
set $D$ of $G_{n,3}$ such that $|D\cap Y_n|=3$. Then $D$ is still a
total dominating set of $G_{n+1,3}$. By Lemma~\ref{lem2.1},
$n=\gamma_t(G_{n,3})=|D|=\gamma_t(G_{n+1,3})=n+1$, a contradiction.
Therefore, $|D\cap Y_n|\leq 2$.
\end{pf}


\begin{lem}\label{lem4.2}
Let $D$ be a total dominating set of $G_{n,3}$. If at least one of
$x_{n1}$ and $x_{n3}$ is in $D$, then $|D|\geq n+1$.
\end{lem}

\begin{pf}
Let $D$ be a total dominating set of $G_{n,3}$.

We first consider that both $x_{n1}$ and $x_{n3}$ are in $D$. If
$X_{n2}\in D$, then $|D\cap Y_n|=3$. By Lemma~\ref{lem4.1}, $D$ is
not a minimum total dominating set of $G_{n,3}$. Thus, by
Lemma~\ref{lem2.1}, we have $|D|\geq \gamma_t(G_{n,3})+1=n+1$.

Assume $x_{n2}\notin D$ below. Since $x_{n2}$ is not in $D$, both
$x_{(n-1)1}$ and $x_{(n-1)3}$ must be in $D$, which dominate
$x_{n1}$ and $x_{n3}$, respectively. Let
$D'=(D\setminus\{x_{n1},x_{n3}\})\cup\{x_{(n-1)2}\}$. Then $D'$ is
still a total dominating set of $G_{n,3}$, and $|D'|<|D|$. By
Lemma~\ref{lem2.1}, $|D|\geq |D'|+1\geq \gamma_t(G_{n,3})+1=n+1$.





We now consider that only one of $x_{n1}$ and $x_{n3}$ is in $D$.
Without loss of generality, we can assume $x_{n1}\in D$ and
$x_{n3}\notin D$. We prove $|D|\geq n+1$ by induction on $n$.

It is clear that $|D|\geq 3$ for $n=2$. Suppose $|D|\geq k+1$ for
any integer $k<n$. We prove that $|D|\geq n+1$ for $n\geq 3$. We can
assume that $|D\cap Y_i|\leq 2$ for each $i=2,3,\ldots, n-1$ since
if $|D\cap Y_i|= 3$ for some $i$ with $2\leq i\leq n-1$, then
$(D\setminus\{x_{i1}, x_{i3}\})\cup \{x_{(i-1)2}, x_{(i+1)2}\}$ is
still a total dominating set of $G_{n,3}$ with the cardinality at
most $|D|$.

If $x_{i2}\notin D$ for each $i=2,3,\ldots,n-1$, then each vertex in
$D\setminus \{x_{n1}\}$ can totally dominate at most three vertices,
and $x_{n1}$ can dominate only two vertices. Thus, $D$ can totally
dominate at most $2+3(|D|-1)$ vertices. On the other hand, $D$ can
totally dominate all $3n$ vertices. From the two facts, we can
deduce $3n\leq 2+3(|D|-1)$, which yields $|D|\geq n+1$.

Now assume $x_{i2}\in D$ for some ${i}$ with $2\leq {i} \leq n-1$.
Let ${i_0}$ be the largest index such that $x_{{i_0}2}\in D$ for
$2\leq {i_0} \leq n-1$. If $x_{n2}\in D$ and ${i_0}=n-1$, then
$D\setminus \{x_{n1}\}$ is still a total dominating set of
$G_{n,3}$. Thus, $|D|\ge \gamma_t(G_{n,3})+1=n+1$ by
Lemma~\ref{lem2.1}. We assume ${i_0}\neq n-1$ if $x_{n2}\in D$ in
the following discussion. There are two cases

\vskip10pt

\begin{description}
\item [Case 1]
$D\cap Y_{i_0}=\{x_{{i_0}2}\}$.

In this case, since $x_{({i_0}+1)2}\notin D$ by the maximality of
${i_0}$, $x_{({i_0}-1)2}\in D$, which dominates the vertex
$x_{{i_0}2}$. Let $D_1=D\cap V(G_{{i_0},3})$. Then $D_1$ is a total
dominating set of $G_{{i_0},3}$.

\begin{center}
\begin{figure}[ht]
\begin{pspicture}(-0.7,0.)(10.5,4.0)
\psset{unit=1.3cm}

\cnode(1,1){2pt}{11}\rput(1,.8){\scriptsize $x_{11}$}
\cnode(1,2){2pt}{12}\rput(1.2,2.2){\scriptsize $x_{12}$}
\cnode(1,3){2pt}{13}\rput(1,3.2){\scriptsize $x_{13}$}

\cnode(2,1){2pt}{21}  \cnode(2,2){2pt}{22}  \cnode(2,3){2pt}{23}

\cnode(3,1){2pt}{31}  \cnode(3,2){2pt}{32}  \cnode(3,3){2pt}{33}

\cnode(4,1){2pt}{41}\rput(4,0.8){\scriptsize $x_{({i_0}-2)1}$}
\cnode(4,2){2pt}{42}\rput(4.5,2.2){\scriptsize $x_{({i_0}-2)2}$}
\cnode(4,3){2pt}{43}\rput(4,3.2){\scriptsize $x_{({i_0}-2)3}$}

\cnode(5,1){2pt}{51}\rput(5,0.8){\scriptsize $x_{({i_0}-1)1}$}
\cnode*(5,2){2pt}{52}\rput(5.5,2.2){\scriptsize $x_{({i_0}-1)2}$}
\cnode(5,3){2pt}{53}\rput(5,3.2){\scriptsize $x_{({i_0}-1)3}$}

\cnode(6,1){2pt}{61}\rput(6,0.8){\scriptsize $x_{{i_0}1}$}
\cnode*(6,2){2pt}{62}\rput(6.3,2.2){\scriptsize $x_{{i_0}2}$}
\cnode(6,3){2pt}{63}\rput(6,3.2){\scriptsize $x_{{i_0}3}$}

\cnode(7,1){2pt}{71}\rput(7,0.8){\scriptsize $x_{({i_0}+1)1}$}
\cnode(7,2){2pt}{72}\rput(7.5,2.2){\scriptsize $x_{({i_0}+1)2}$}
\cnode(7,3){2pt}{73}\rput(7,3.2){\scriptsize $x_{({i_0}+1)3}$}

\cnode(8,1){2pt}{81}  \cnode(8,2){2pt}{82}   \cnode(8,3){2pt}{83}

\cnode(9,1){2pt}{91}  \cnode(9,2){2pt}{92}   \cnode(9,3){2pt}{93}

\cnode*(10,1){2pt}{101}\rput(10,0.8){\scriptsize $x_{n1}$}
\cnode(10,2){2pt}{102}\rput(10.3,2.2){\scriptsize $x_{n2}$}
\cnode(10,3){2pt}{103}\rput(10,3.2){\scriptsize $x_{n3}$}

\ncline{11}{12} \ncline{12}{13}  \ncline{11}{21} \ncline{12}{22}
\ncline{13}{23} \ncline{21}{22} \ncline{22}{23}
\ncline[linestyle=dotted]{21}{31} \ncline[linestyle=dotted]{22}{32}
\ncline[linestyle=dotted]{23}{33} \ncline{31}{32} \ncline{32}{33}
\ncline{31}{41} \ncline{32}{42}  \ncline{33}{43} \ncline{41}{42}
\ncline{42}{43}  \ncline{41}{51} \ncline{42}{52}  \ncline{43}{53}
\ncline{51}{52} \ncline{52}{53}  \ncline{51}{61} \ncline{52}{62}
\ncline{53}{63} \ncline{61}{62} \ncline{62}{63}  \ncline{61}{71}
\ncline{62}{72}  \ncline{63}{73} \ncline{71}{72} \ncline{72}{73}
\ncline{71}{81} \ncline{72}{82}  \ncline{73}{83} \ncline{81}{82}
\ncline{82}{83}  \ncline[linestyle=dotted]{81}{91}
\ncline[linestyle=dotted]{82}{92}  \ncline[linestyle=dotted]{83}{93}
\ncline{91}{92} \ncline{92}{93}  \ncline{91}{101} \ncline{92}{102}
\ncline{93}{103} \ncline{101}{102}  \ncline{102}{103}

\rput[b]{-90}(.8,0.4){$\left.\rule{0mm}{35mm}\right\}~\rput{90}{G_{{i_0},3}}$}
\rput[b]{-90}(4.9,0.2){$\left.\rule{0mm}{35mm}\right\}~\rput{90}{H_{n-({i_0}-2),3}}$}
\end{pspicture}
\caption{\label{f3}\footnotesize  Two subgraphs $G_{{i_0},3}$ and
$H_{n-({i_0}-2),3}$ of $G_{n,3}$}
\end{figure}
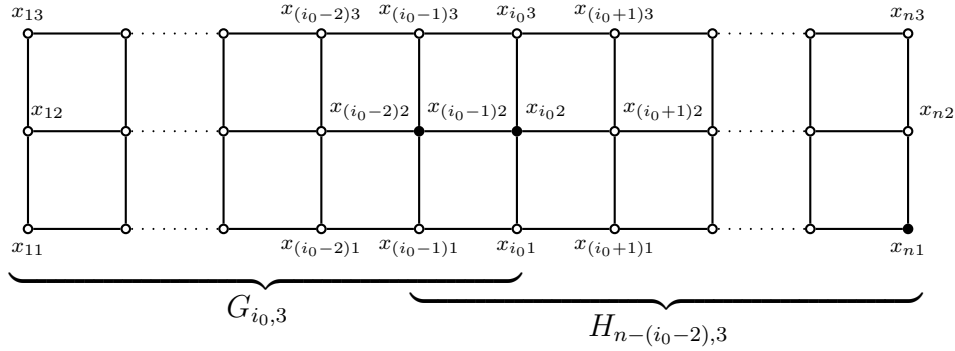
\end{center}

If $D\cap Y_{{i_0}-1}=\{x_{({i_0}-1)2}\}$, then $D_2=D\cap
V(H_{n-({i_0}-2),3})$ is a total dominating set of
$H_{n-({i_0}-2),3}$, and $|D_1\cap D_2|=|\{x_{{i_0-1}2},
x_{i_02}\}|=2$ (see Figure~\ref{f3}). Since $H_{n-({i_0}-2),3}\cong
G_{n-({i_0}-2),3}$ and $D_2$ satisfies the condition in the lemma
(i.e., $x_{n1}\in D_2$), by the induction hypothesis $|D_2|\geq
(n-{i_0}+2)+1=n-{i_0}+3$. By Lemma~\ref{lem2.1}, $|D_1|\geq {i_0}$.
Thus,
%
 $$
 |D|\geq |D_1|+|D_2|-2\geq {i_0}+n-{i_0}+3-2=n+1.
 $$

If $|D\cap Y_{{i_0}-1}|=2$, then $D_3=D\cap V(H_{n-{(i_0-2)},3})$ is
a total dominating set of $H_{n-({i_0}-2),3}$. Let $D_4=D\cap
V(G_{i_0-1},3)$. Then $|D_4\cap D_3|=2$. Since both of $D_3$ and
$D_4$ satisfy the condition in the lemma, by the induction
hypothesis, $|D_4|\geq {i_0}$ and $|D_3|\geq n-({i_0}-2)+1$. Thus,
 $$
 |D|\geq |D_4|+|D_3|-2\geq {i_0}+n-{i_0}+3-2 =n+1.
 $$

\item [Case 2]
If $|D\cap Y_{i_0}|=2$, then $D_5=D\cap V(G_{{i_0},3})$ is a total
dominating set of $G_{{i_0},3}$, $D_6=D\cap V(H_{n-({i_0}-1),3})$ is
a total dominating set of $H_{n-({i_0}-1),3}$ and $|D_5\cap D_6|=2$.
Since both $D_4$ and $D_5$ satisfy the condition in the lemma, by
the induction hypothesis, $|D_5|\geq {i_0}+1$ and $|D_6|\geq
n-{i_0}+1+1=n-{i_0}+2$. Thus,
 $$
 |D|\geq |D_5|+|D_6|-2\geq {i_0}+1+n-{i_0}+2-2=n+1.
 $$
\end{description}

The proof of the lemma is complete.
\end{pf}


\begin{thm}\label{thm4.1}
$b_{t}(G_{n,3})=1$ for $n\geq 2$.
\end{thm}

\begin{pf}
Let $H=G_{n,3}-x_{(n-1)2}x_{n2}$ and $D$ be a minimum total
dominating set of $H$. Whether $x_{n2}$ is in $D$  or not, at least
one of two vertices $x_{n1}$ and $x_{n3}$ is in $D$, which dominates
$x_{n2}$ in $H$. By Lemma~\ref{lem4.2}, $|D|\geq n+1$. Combining
this fact with Lemma~\ref{lem2.1}, we have $\gamma_t(H)=|D|\geq
n+1=\gamma_t(G_{n,3})+1$. Therefore, $b_{t}(G_{n,3})=1$.
\end{pf}


\section{The total bondage number of $G_{n,4}$}

In this section, we determine the exact value of $b_t(G_{n,4})$ for
$n\equiv 1,\,4\,({\rm mod}\, 5)$, and establish the upper bounds of
$b_t(G_{n,4})$  for $n\equiv 0,2,3\,({\rm mod}\, 5)$.

\begin{lem}\label{lem5.1}
$b_{t}(G_{n,4})=1$ for $n\geq 7$ and $n\equiv 1\,({\rm mod}\, 5)$.
\end{lem}

\begin{pf}
Let
$D$ be a minimum total dominating set of $G_{n,4}-x_{n2}x_{n3}$. It
is easy to see that $|D\cap (Y_{n-1}\cup Y_n)|\geq 4$. Thus,
$|D|\geq 4+|D\cap V(G_{n-2,4})|$. When $n\ge 7$, $n-3\equiv 3\,({\rm
mod}\, 5)$ and $n-3\ge 4$. By Lemma~\ref{lem2.1}, we have
$\gamma_{t}(G_{n-3,4})=\lfloor\frac{6(n-3)+8}{5}\rfloor+1$. Thus, by
Lemma \ref{lem2.2}, we have
 $$
  \begin{array}{rl}
 |D\cap V(G_{n-2,4})|&\geq
 \gamma_{t}(G_{n-3,4})\\
 &=\left\lfloor\frac{6(n-3)+8}{5}\right\rfloor+1=\left\lfloor\frac{6n+8}{5}\right\rfloor-2
 \end{array}
 $$
and, hence,
 $$
 \gamma_{t}(G_{n,4}-x_{n2}x_{n3})=|D|\geq 2+\left\lfloor\frac{6n+8}{5}\right\rfloor>\gamma_{t}(G_{n,4}).
 $$
Therefore, $b_{t}(G_{n,4})=1$.
\end{pf}


\vskip6pt

To determine $b_{t}(G_{n,4})$ for $n\equiv 4\,({\rm mod}\,5)$, we
state two simple observations, see Figure~\ref{f4} for $n=9$.


\begin{center}
\begin{figure}[ht]
\begin{pspicture}(-8,0.3)(9,3.5)
\psset{unit=.8cm}

\cnode*(-1,1){2pt}{l11}\rput(-1,.7){\scriptsize $x_{91}$}
\cnode(-1,2){2pt}{l12}  \cnode(-1,3){2pt}{l13}
\cnode*(-1,4){2pt}{l14}\rput(-1,4.3){\scriptsize  $x_{94}$}

\cnode*(-2,1){2pt}{l21}\rput(-2,.7){\scriptsize $x_{81}$}
\cnode(-2,2){2pt}{l22}   \cnode(-2,3){2pt}{l23}
\cnode*(-2,4){2pt}{l24}\rput(-2,4.3){\scriptsize $x_{84}$}

\cnode(-3,1){2pt}{l31}\rput(-3,.7){\scriptsize $x_{71}$}
\cnode(-3,2){2pt}{l32}  \cnode(-3,3){2pt}{l33}
\cnode(-3,4){2pt}{l34}\rput(-3,4.3){\scriptsize $x_{74}$}

\cnode(-4,1){2pt}{l41}\rput(-4,.7){\scriptsize $x_{61}$}
\cnode*(-4,2){2pt}{l42}  \cnode*(-4,3){2pt}{l43}
\cnode(-4,4){2pt}{l44}\rput(-4,4.3){\scriptsize $x_{64}$}

\cnode(-5,1){2pt}{l51}\rput(-5,.7){\scriptsize $x_{51}$}
\cnode(-5,2){2pt}{l52}  \cnode(-5,3){2pt}{l53}
\cnode(-5,4){2pt}{l54}\rput(-5,4.3){\scriptsize $x_{54}$}

\cnode*(-6,1){2pt}{l61}\rput(-6,.7){\scriptsize $x_{41}$}
\cnode(-6,2){2pt}{l62}  \cnode(-6,3){2pt}{l63}
\cnode*(-6,4){2pt}{l64}\rput(-6,4.3){\scriptsize $x_{44}$}

\cnode*(-7,1){2pt}{l71}\rput(-7,.7){\scriptsize $x_{31}$}
\cnode(-7,2){2pt}{l72}  \cnode(-7,3){2pt}{l73}
\cnode*(-7,4){2pt}{l74}\rput(-7,4.3){\scriptsize $x_{34}$}

\cnode(-8,1){2pt}{l81}\rput(-8,.7){\scriptsize $x_{21}$}
\cnode(-8,2){2pt}{l82}  \cnode(-8,3){2pt}{l83}
\cnode(-8,4){2pt}{l84}\rput(-8,4.3){\scriptsize $x_{24}$}

\cnode(-9,1){2pt}{l91}\rput(-9,.7){\scriptsize $x_{11}$}
\cnode*(-9,2){2pt}{l92}  \cnode*(-9,3){2pt}{l93}
\cnode(-9,4){2pt}{l94}\rput(-9,4.3){\scriptsize $x_{14}$}

\ncline{l11}{l12} \ncline{l12}{l13}  \ncline{l11}{l21}
\ncline{l12}{l22}  \ncline{l13}{l23}  \ncline{l13}{l14}
\ncline{l14}{l24} \ncline{l21}{l22} \ncline{l22}{l23}
\ncline{l21}{l31} \ncline{l22}{l32}  \ncline{l23}{l33}
\ncline{l23}{l24}  \ncline{l24}{l34} \ncline{l31}{l32}
\ncline{l32}{l33}  \ncline{l31}{l41} \ncline{l32}{l42}
\ncline{l33}{l43}  \ncline{l33}{l34}  \ncline{l34}{l44}
\ncline{l41}{l42} \ncline{l42}{l43}  \ncline{l41}{l51}
\ncline{l42}{l52}  \ncline{l43}{l53}  \ncline{l43}{l44}
\ncline{l44}{l54} \ncline{l51}{l52} \ncline{l52}{l53}
\ncline{l51}{l61} \ncline{l52}{l62}  \ncline{l53}{l63}
\ncline{l53}{l54}  \ncline{l54}{l64} \ncline{l61}{l62}
\ncline{l62}{l63}  \ncline{l61}{l71} \ncline{l62}{l72}
\ncline{l63}{l73}  \ncline{l63}{l64}  \ncline{l64}{l74}
\ncline{l71}{l72} \ncline{l72}{l73}  \ncline{l71}{l81}
\ncline{l72}{l82}  \ncline{l73}{l83}  \ncline{l73}{l74}
\ncline{l74}{l84} \ncline{l81}{l82} \ncline{l82}{l83}
\ncline{l81}{l91} \ncline{l82}{l92}  \ncline{l83}{l93}
\ncline{l83}{l84}  \ncline{l84}{l94} \ncline{l91}{l92}
\ncline{l92}{l93}  \ncline{l93}{l94}

\cnode*(1,1){2pt}{11}\rput(1,.7){\scriptsize $x_{11}$}
\cnode*(1,2){2pt}{12}  \cnode(1,3){2pt}{13}
\cnode(1,4){2pt}{14}\rput(1,4.3){\scriptsize $x_{14}$}

\cnode(2,1){2pt}{21}\rput(2,.7){\scriptsize $x_{21}$}
\cnode(2,2){2pt}{22}  \cnode(2,3){2pt}{23}
\cnode*(2,4){2pt}{24}\rput(2,4.3){\scriptsize $x_{24}$}

\cnode(3,1){2pt}{31}\rput(3,.7){\scriptsize $x_{31}$}
\cnode(3,2){2pt}{32}  \cnode(3,3){2pt}{33}
\cnode*(3,4){2pt}{34}\rput(3,4.3){\scriptsize $x_{34}$}

\cnode*(4,1){2pt}{41}\rput(4,.7){\scriptsize $x_{41}$}
\cnode*(4,2){2pt}{42}  \cnode(4,3){2pt}{43}
\cnode(4,4){2pt}{44}\rput(4,4.3){\scriptsize $x_{44}$}

\cnode(5,1){2pt}{51}\rput(5,.7){\scriptsize $x_{51}$}
\cnode(5,2){2pt}{52}  \cnode(5,3){2pt}{53}
\cnode(5,4){2pt}{54}\rput(5,4.3){\scriptsize $x_{54}$}

\cnode(6,1){2pt}{61}\rput(6,.7){\scriptsize $x_{61}$}
\cnode(6,2){2pt}{62}  \cnode*(6,3){2pt}{63}
\cnode*(6,4){2pt}{64}\rput(6,4.3){\scriptsize $x_{64}$}

\cnode*(7,1){2pt}{71}\rput(7,.7){\scriptsize $x_{71}$}
\cnode(7,2){2pt}{72}  \cnode(7,3){2pt}{73}
\cnode(7,4){2pt}{74}\rput(7,4.3){\scriptsize $x_{74}$}

\cnode*(8,1){2pt}{81}\rput(8,.7){\scriptsize $x_{81}$}
\cnode(8,2){2pt}{82}  \cnode(8,3){2pt}{83}
\cnode(8,4){2pt}{84}\rput(8,4.3){\scriptsize $x_{84}$}

\cnode(9,1){2pt}{91}\rput(9,.7){\scriptsize $x_{91}$}
\cnode(9,2){2pt}{92}  \cnode*(9,3){2pt}{93}
\cnode*(9,4){2pt}{94}\rput(9,4.3){\scriptsize $x_{94}$}

\ncline{11}{12} \ncline{12}{13}  \ncline{11}{21} \ncline{12}{22}
\ncline{13}{23}  \ncline{13}{14}  \ncline{14}{24} \ncline{21}{22}
\ncline{22}{23}  \ncline{21}{31} \ncline{22}{32}  \ncline{23}{33}
\ncline{23}{24}  \ncline{24}{34} \ncline{31}{32} \ncline{32}{33}
\ncline{31}{41} \ncline{32}{42}  \ncline{33}{43}  \ncline{33}{34}
\ncline{34}{44} \ncline{41}{42} \ncline{42}{43}  \ncline{41}{51}
\ncline{42}{52}  \ncline{43}{53}  \ncline{43}{44}  \ncline{44}{54}
\ncline{51}{52} \ncline{52}{53}  \ncline{51}{61} \ncline{52}{62}
\ncline{53}{63}  \ncline{53}{54}  \ncline{54}{64} \ncline{61}{62}
\ncline{62}{63}  \ncline{61}{71} \ncline{62}{72}  \ncline{63}{73}
\ncline{63}{64}  \ncline{64}{74} \ncline{71}{72} \ncline{72}{73}
\ncline{71}{81} \ncline{72}{82}  \ncline{73}{83}  \ncline{73}{74}
\ncline{74}{84} \ncline{81}{82} \ncline{82}{83}  \ncline{81}{91}
\ncline{82}{92}  \ncline{83}{93}  \ncline{83}{84}  \ncline{84}{94}
\ncline{91}{92} \ncline{92}{93}  \ncline{93}{94}

\end{pspicture}
\caption{\label{f4}\footnotesize Two minimum total dominating sets
(bold vertices) of $G_{9,4}$ defined in Proposition~\ref{pps5.1} and
Proposition~\ref{pps5.2}, respectively.}
\end{figure}
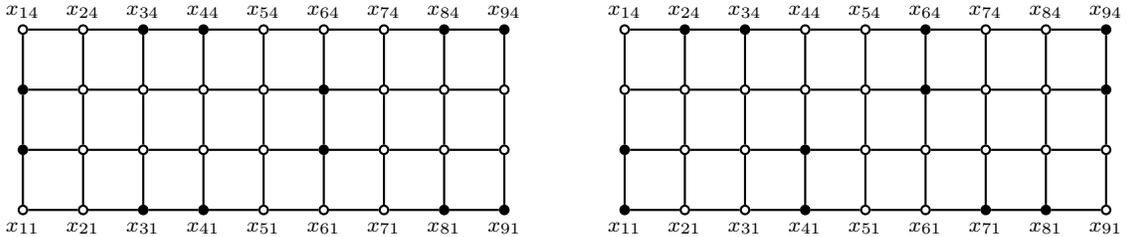
\end{center}

\begin{pps}\label{pps5.1}
For $n\equiv 4\,({\rm mod}\,5)$, both\\
 $
D=\{x_{i2},x_{i3},x_{(i+2)1},x_{(i+3)1},x_{(i+2)4},x_{(i+3)4}:i\equiv
1\,({\rm mod}\, 5), 1\le i\le n-3\}
 $
and\\
 $
D'=\{x_{i2},x_{i3},x_{(i-3)1},x_{(i-2)1},x_{(i-3)4},x_{((i-2))4}:i\equiv
4\,({\rm mod}\,5), 4\le i\le n\}
 $\\
are minimum total dominating sets of $G_{n,4}$.

\end{pps}

\begin{pps}\label{pps5.2}
For $n\equiv 4\,({\rm mod}\,5)$, both\\
$D=\{x_{i1},x_{i2},x_{(i+3)1},x_{(i+3)2},x_{(i+1)4},x_{(i+2)4}: i\equiv 1\,({\rm mod}\,10), 1\le i \le n-3\} \\
~~~~~~~~\cup
\{x_{j3},x_{j4},x_{(j+3)3},x_{(j+3)4},x_{(j+1)1},x_{(j+2)1}: j\equiv
6\,({\rm mod}\,10), 1\le j \le n-3\}$ and\\
$D'=\{x_{i3},x_{i4},x_{(i+3)3},x_{(i+3)4},x_{(i+1)1},x_{(i+2)1}: i\equiv 1\,({\rm mod}\,10), 1\le i \le n-3\} \\
~~~~~~~~\cup \{x_{j1},x_{j2},x_{(j+3)1},x_{(j+3)2},x_{(j+1)4},x_{(j+2)4}: j\equiv 6\,({\rm mod}\,10), 1\le j \le n-3\}$\\
are minimum total dominating sets of $G_{n,4}$.
\end{pps}

\begin{lem}\label{lem5.2}
$b_{t}(G_{n,4})=2$ for $n\equiv 4\,({\rm mod}\,5)$.
\end{lem}

\begin{pf}
We can direct check $b_{t}(G_{4,4})=2$, and assume $n\geq 9$ below.
For any edge $e\in E(G_{n,4})$, it is easy to verify that $D$ or
$D'$ defined in Proposition~\ref{pps5.1} or Proposition~\ref{pps5.2}
is also a minimum total dominating set of $G_{n,4}-e$. Thus,
$b_{t}(G_{n,4})\ge 2$. We now prove that $b_{t}(G_{n,4})\le 2$.

Let $H=G_{n,4}-x_{(n-1)1}x_{n1}-x_{(n-1)2}x_{n2}$ and let $S$ be a
minimum total dominating set of $H$. Then the vertex  $x_{n2}$ must
be in $S$ to dominate $x_{n1}$, and at least one of $x_{n1}$ and
$x_{n3}$ must be in $S$ to dominate $x_{n2}$ in $H$, that is,
$|Y_n\cap S|\ge 2$.

If $|Y_n\cap S|\ge 3$ then $|S|\geq |S\cap V(G_{n-1,4})|+3$. By
Lemma \ref{lem2.2} and Lemma~\ref{lem2.1}, we have
 $$
  \begin{array}{rl}
 |S\cap V(G_{n-1,4})|&\geq \gamma_{t}(G_{n-2,4})\\
 &= \left\lfloor\frac{6(n-2)+8}{5}\right\rfloor=
 \left\lfloor\frac{6n+8}{5}\right\rfloor-2
 \end{array}
 $$
and, hence,
 $$
 \gamma_{t}(H)=|S|\geq \gamma_{t}(G_{n-2,4})+3
 \geq 1+\left\lfloor\frac{6n+8}{5}\right\rfloor>\gamma_{t}(G_{n,4}).
 $$

If $|Y_n\cap S|=2$, then $Y_n\cap S$ can totally dominate at most
one vertex in $G_{n-1,4}$, that is, $x_{(n-1)3}$ if so. Thus,
$(S\cap V(G_{n-1,4}))\cup \{x_{(n-1)2}\}$ is a total dominating set
of $G_{n-1,4}$, which implies $|S|\geq |(S\cap V(G_{n-1,4}))\cup
\{x_{(n-1)2}\}|+1\geq \gamma_{t}(G_{n-1,4})+1$. By
Lemma~\ref{lem2.1}, we have
 $$
 \gamma_{t}(H)=|S|\ge \gamma_{t}(G_{n-1,4})+1>\gamma_{t}(G_{n,4}).
 $$

Therefore, $b_{t}(G_{n,4})\leq 2$.
\end{pf}


\begin{lem}\label{lem5.3}
$b_{t}(G_{n,4})\leq 3$ for $n\equiv 2\,({\rm mod}\, 5)$.
\end{lem}

\begin{pf}
$b_{t}(G_{4,2})=3$ by Theorem~\ref{thm3.1}, and $b_{t}(G_{4,7})\leq
3$ by checking direct. Assume $n\geq 12$ below. Let
$H=G_{n,4}-x_{(n-1)1}x_{n1}-x_{(n-1)2}x_{n2}-x_{n2}x_{n3}$ and let
$S$ be a minimum total dominating set of $H$. Since $x_{n1}x_{n2}$
is an isolated edge in $H$, both $x_{n1}$ and $x_{n2}$ must be in
$S$. To dominate the three vertices $x_{(n-1)1}$, $x_{n3}$ and
$x_{n4}$, we need at least three other vertices in $S$. In other
words, $|(Y_{n-2}\cup Y_{n-1}\cup Y_n)\cap S|\geq 5$. Thus,
 \begin{equation}\label{e5.1}
 |S|\geq 5+|S\cap V(G_{n-3,4})|.
 \end{equation}
By Lemma \ref{lem2.2} and Lemma \ref{lem2.1}, when $n\geq 8$ and
$n\equiv 2\,({\rm mod}\, 5)$, we have
 $$
 \begin{array}{rl}
 |S\cap V(G_{n-3,4})|&\geq  \gamma_{t}(G_{n-4,4})\\
 &=\left\lfloor\frac{6(n-4)+8}{5}\right\rfloor+1=\left\lfloor\frac{6n+8}{5}\right\rfloor-4,
  \end{array}
  $$
that is,
\begin{equation}\label{e5.2}
 |S\cap V(G_{n-3,4})|\geq \left\lfloor\frac{6n+8}{5}\right\rfloor-4.
 \end{equation}
It follows from (\ref{e5.1}), (\ref{e5.2}) and Lemma~\ref{lem2.1}
that
 $$
  \gamma_{t}(H)=|S|\geq 1+\left\lfloor\frac{6n+8}{5}\right\rfloor\\
  >\gamma_{t}(G_{n,4}).
 $$

Therefore, $b_{t}(G_{n,4})\leq 3$.
\end{pf}


\begin{lem}\label{lem5.4}
$b_{t}(G_{n,4})\leq 4$ for $n\equiv 0,3\,({\rm mod}\, 5)$.
\end{lem}

\begin{pf}
By direct checking we have $b_{t}(G_{5,4}))\leq 4$ and
$b_{t}(G_{8,4}))\leq 4$. Assume $n\geq 10$ below. Let
$H=G_{n,4}-x_{(n-5)1}x_{(n-4)1}-x_{(n-5)2}x_{(n-4)2}-x_{(n-5)3}x_{(n-4)3}-x_{(n-5)4}x_{(n-4)4}$
and let $S$ be a minimum total dominating set of $H$. Then $H$
consists of two grid subgraphs $G_{n-5,4}$ and $H_{5,4}$, and so
\begin{equation}\label{e5.3}
|S|=\gamma_{t}(G_{n-5,4})+ \gamma_{t}(H_{5,4}).
 \end{equation}
Since $n-5\equiv 0,3\,({\rm mod}\, 5)$ and $5\equiv 0\,({\rm mod}\,
5)$, by Lemma~\ref{lem2.1}, we have
 \begin{equation}\label{e5.4}
 \begin{array}{rl}
 \gamma_{t}(G_{n-5,4})+ \gamma_{t}(H_{5,4})
   &=\left\lfloor\frac{6(n-5)+8}{5}\right\rfloor+1+8\\
   &=\left\lfloor\frac{6n+8}{5}\right\rfloor+3
   >\gamma_{t}(G_{n,4}).
   \end{array}
 \end{equation}
Combining (\ref{e5.3}) with (\ref{e5.4}), we have
 $$
 \gamma_{t}(H)=\gamma_{t}(G_{n-5,4})+ \gamma_{t}(H_{5,4})>\gamma_{t}(G_{n,4}).
 $$

Therefore, $b_{t}(G_{n,4})\leq 4$.
\end{pf}


\vskip6pt

Summing up the above lemmas, we can state our result, in this
section, as follows.

\begin{thm}\label{thm5.1} For any integer $n\geq 1$, $b_{t}(G_{6,4})=2$, and
 $$
 b_{t}(G_{n,4})\left\{ \begin{array}{ll}
 =1\ & {\rm if}\ n\equiv 1\,({\rm mod}\, 5)\ {\rm and}\ n\ne 6;\\
 =2\ & {\rm if}\  n\equiv 4\,({\rm mod}\, 5);\\
 \leq 3\ & {\rm if}\ n\equiv 2\,({\rm mod}\, 5);\\
 \leq 4\ & {\rm if}\ n\equiv 0,3\,({\rm mod}\, 5).\\
 \end{array} \right.
 $$
\end{thm}

\section{Concluding Remarks}

In this paper, we investigate the total bondage number
$b_t(G_{n,m})$ of an $(n,m)$-grid graph $G_{n,m}$ for $2\leq m\leq
4$, completely determine the exact values of $b_t(G_{n,2})$ and
$b_t(G_{n,3})$. We also partially determine the exact values of
$b_t(G_{n,4})$, and establishes the upper bounds of $b_t(G_{n,4})$
for otherwise. We have attempted to decrease the two upper bounds
given in Theorem~\ref{thm5.1} for $n\equiv 2\,({\rm mod}\, 5)$ and
$n\equiv 0,3\,({\rm mod}\, 5)$ or to prove that they are tight when
$n$ is large enough, but we have not been able to bring home the
bacon.
Noting the two upper bounds are tight for some small $n$'s, we guess
that the two upper bounds are tight for $n\geq 7$. To prove this
conjecture, it may be necessary to find a new method since,
according to our way, the removal of any three edges results in many
complicated cases. We also have tried to discuss $b_t(G_{n,m})$ for
general $n$ and $m$, but it strongly depends on the value of
$\gamma_t(G_{n,m})$, which has not been determined as yet. Thus, it
may also be necessary to determine the value of $\gamma_t(G_{n,m})$
for general $n$ and $m$. These questions are our further work.

\section*{Acknowledgements}

The authors would like to thank the anonymous referees for their
helpful comments and kind suggestions on the original manuscript,
which resulted in this final version.

\end{document}